%%%%%%%%%%%%%%%%%%%%%%%%%%%%%%%%%%%%%%%%%%%%%%%%%%%%%%%%%%%%%%%%
%%%%
%%%%           AMSTeX
%%%%
%%%%%%%%%%%%%%%%%%%%%%%%%%%%%%%%%%%%%%%%%%%%%%%%%%%%%%%%%%%%%%%%
%%%%%%%%%%%%%%%%%%%%%%%%%%%%%%%%%%%%%%%%%%%%%%%%%%%%%%%%%%%%%%%%
 \input amstex
\documentstyle{amsppt}
\magnification=\magstep1
\parindent=1em
\CenteredTagsOnSplits
\NoBlackBoxes
\nopagenumbers
\NoRunningHeads
\pageno=1
\footline={\hss\tenrm\folio\hss}

        \topmatter
        \title {
     On linear operators with ${\ssize\bold s}$-nuclear adjoints: $0<{\ssize s}\le 1$
               }
        \endtitle
          \author { O.I. Reinov%${{ }^\dag}$
          }  \endauthor
\address\newline
Oleg I. Reinov \newline
%Department of Mathematics\newline
Saint Petersburg  State University\newline
198504 St. Petersburg, Petrodvorets, \newline
28 Universitetskii pr., Russia
\endaddress

\email
orein51\@mail.ru
\endemail

%\thanks
%$%{{ }^\dag}$This work was done with partial support of the Ministry of the
%general and professional education of Russia (Grant 97-0-1.7-36) and
%FCP ``Integracija", reg. No. 326.53.
%\endthanks

\abstract
If $ s\in (0,1]$ and $ T$ is a linear operator with $ s$-nuclear
adjoint from a Banach space $ X$ to a Banach space $ Y$ and if one
of the spaces $ X^*$ or $ Y^{***}$ has the approximation property of order $s,$ \, $AP_s,$ then
the operator $ T$ is nuclear.
The result is in a sense exact. For example, it is shown that for each $r\in (2/3, 1]$
there exist a Banach space $Z_0$ and a non-nuclear operator
$ T: Z_0^{**}\to Z_0$ so that
 $ Z_0^{**}$ has
a Schauder basis,  $ Z_0^{***}$ has the $AP_s$ for every $s\in (0,r)$
and $T^*$ is $r$-nuclear.
\endabstract
    \endtopmatter

\document
\baselineskip=18pt

\footnote""{${ }^\ddag$
}
\footnote""{${ }$
}

%\def\al{\alpha}
%\def\ot{\otimes}
%\def\wh{\widehat}
%\def\ffi{\varphi}
%\def\wt{\widetilde}
%\def\small{\smallpagebreak}
%\def\sbs{\subset}
%\def\({\left(}
%\def\){\right)}
%\def\({\left(}
%\def\){\right)}
%\def\[{\left[}
%\def\]{\right]}
%\def\tr{\operatorname{trace}\,}

%\def\<{\langle}
  %   \def\>{\rangle}
%\def\id{\operatorname{id}}
%\def\e{\varepsilon}

%%%%%%%%%%%%%%%%%%%%%%%%%%%%%%%%%%%%%%%%%%%%%%%

\centerline{\bf \S1. \, Introduction}
\smallpagebreak

We will be interested in the following question from [6, Problem 10.1]:           %R!!!
suppose $T$ is a (bounded linear) operator acting between Banach spaces $X$ and $Y,$
and let $s\in(0,1).$ Is it true that if $T^*$ is $s$-nuclear then $T$ is $s$-nuclear too?

As well known, for $s=1$ the negative answer was obtained for the first time
by Figiel and Johnson in [4].                   %R!!!
For $s\in (2/3,1]$ the negative answer can be found, e.g., in [14].      %R!!!
Here, we are going to give some (partially) positive results in this direction as well as to show the sharpness of them.

It is not difficult to see that if $T^*\in N_s(Y^*,X^*)$ then $T\in N_p(X,Y)$ with $1/s=1/p+1/2$
(what is, surely, must be known; see, e.g., [13], [14]). This is the best possible result in the scale
of $q$-nuclear operators:                     %RR!!!
the sharpness of the last assertion, for $s\in (2/3, 1],$
can be found, for instance, in [14].                  %R!!!

Below, we consider a little bit different question:
under which assumptions on Banach spaces $X$ and $Y$

$(*)$\  an operator $T\in L(X,Y)$ is nuclear if its adjoint $T^*$ is $s$-nuclear?

One of the possibilities for getting some positive answers to $(*)$
is to apply the so-called approximation properties (the $AP_s)$
of order $s,$\, $s\in(2/3,1]$
(we assume that $s>2/3$ since for $s\le 2/3$ the answer is evident).

We will prove that $(*)$ is true if either $X^*$ or $Y^{***}$ has the $AP_s$ (Theorem 1).
Some examples, given after Theorem 1, will show that those assumptions are, in a sense,
necessary (for example, it is not enough to assume that either $X$ or $Y^{**},$ or both of the spaces
have the approximation properties, even in the sense of A.Grothendieck).
Let us note that the case where $s=1$ was firstly investigated 
in the paper  [8] by Eve Oja and Oleg Reinov.

%%%%%%%%%%%%%%%%%%%%%%%%   Perehod k Temam
\vskip0.5cm

%{\bf Notations, preliminaries.}\

\centerline{\bf \S2. \, Notations and  preliminaries}
\smallpagebreak

For the references of different things, one can see [9].  %R!!!
All the spaces under considerations $(X,Y,\dots)$
are Banach, all linear mappings (operators) are continuous; as usual, $X^*, X^{**}, \dots$
are Banach duals (to $X$), and $x', x'', \dots$ (or $y', \dots)$ are the functionals on $X, X^*,\dots$
(or on $Y,\dots).$ By $\pi_Y$ we denote the natural isometric injection of $Y$ into its second dual.
If $x\in X, x'\in X^*$ then $\langle x,x'\rangle=\langle x',x\rangle=x'(x).$
$L(X,Y)$ stands for the Banach space of all linear bounded operators from $X$ to $Y.$

Recall that an operator $T:X\to Y$ is $s$-nuclear $(0<s\le1)$ if it is of the form
$$
 Tx=\sum_{k=1}^\infty \langle x'_k,x\rangle y_k
$$
for $x\in X,$ where $(x'_k)\subset X^*, (y_k)\subset Y,\, \sum_k ||x'_k||^s\,||y_k||^s<\infty.$ We use
the notation $N_s(X,Y)$ for the space of all such operators.
Every $s$-nuclear operator is a canonical image of an element of a projective tensor product, namely:
denote by $X^*\widehat\otimes_s Y$ the $s$-projective tensor product of $X^*$ and $Y$ that is
a subspace of the Grothendieck projective tensor product
$X^*\widehat\otimes Y$ $(= X^*\widehat\otimes_1 Y)$\, consisting of all tensor elements $z$ which admit
representations of kind
$$
 z=\sum_{k=1}^\infty x'_k\otimes y_k  \ \text{ with }\ \sum_{k=1}^\infty ||x'_k||^s\,||y_k||^s<\infty.
$$
Then every $s$-nuclear operator from $X$ to $Y$ is an image of an element of $X^*\widehat\otimes_s Y$
via the canonical mappings
$$
 X^*\widehat\otimes_s Y\overset{j_s}\to\to X^*\widehat\otimes_1 Y\overset{j}\to\to L(X,Y).
$$
If $z\in X^*\widehat\otimes Y$ then we denote the corresponding operator (from $X$ to $Y)$
by $\tilde z.$

 We say, following, e.g., [13] or  [14],        %R!!!
 that a Banach space $Y$ has the $AP_s$ (the approximation property of order $s),$
 if for every Banach space $X$ the natural map $jj_s$ (from above) is one-to-one
 (note that $AP_1=AP$ of A. Grothendieck).
 It is the same as to say that the natural map $Y^*\widehat\otimes_s Y \to L(Y,Y)$ is one-to-one.
 Therefore, if $Y$ has the $AP_s,$ we can write $X^*\widehat\otimes_s Y = N_s(X,Y)$
 whichever $X$ was.
    Clearly, $AP_s\implies AP_t$ for $0<t<s\le1.$ Every Banach space has the $AP_{2/3}$
  (Grothendieck's Theorem [5], see also [11], [13] or [14]).         %RR!!!

  Recall that the dual space to $X^*\widehat\otimes_1 Y$ is $L(X,Y^{**})$ and the duality is defined by "trace":
  if $z\in X^*\widehat\otimes_1 Y$ and $U\in L(Y, X^{**}),$ then
  $$
   \langle U,z\rangle := \operatorname{trace}\, U\circ z
  $$
  $(= \sum_k \langle x'_k, Uy_k\rangle$ for a projective representation $z= \sum_k x'_k\otimes y_k).$
  So, the element $z\in X^*\widehat\otimes_s Y$ is zero iff it it zero in the projective tensor product $X^*\widehat\otimes Y$
  iff for every $U\in L(Y, X^{**})$ \, $\operatorname{trace}\, U\circ z=0.$ If $z\in X^*\widehat\otimes Y$ then the corresponding
  operator $\tilde z$ is zero iff for every $R\in Y^*\otimes X$ we have: $\operatorname{trace}\, R\circ z=0$ (evidently).

  Examples of Banach spaces with $AP_s:$ for $s\in[2/3,1]$ and $1/p+1/2=1/s,$ every quotient of any subspace
  of any $L_p$-space (and every subspace of any quotient of any $L_{p'}$-space) has the $AP_s$
  (as well as all their duals; see, e.g., [13] or [14]; for a more general fact, see Lemma 3 below).
  Here $1/p+1/p'=1.$
  All the Banach spaces have $AP_s$ for $s\in (0,2/3];$ but if $2/3\le s_1< s_2\le1$
  then $AP_{s_2}\implies AP_{s_1}$ but $AP_{s_1}$ does not imply $AP_{s_2}.$
  It is known that for every $p\neq 2,$ $p\in [1,\infty],$
  there exists a subspace (a lot of ones) of $l_p$ without the Grothendieck approximation property;
  thus, for example, $AP_1\neq AP_s$ if $s\in (0,1).$ Something (about "$AP_{s_1}\neq AP_{s_2}$")
  will be possible to get from Examples or Theorems 2-3 below.

We will use later the following facts (surely, well known, but maybe not mentioned in the literature):
\smallpagebreak

{\bf Lemma 1.}\
If $T\in L(X, Y^{**})$ then
$||T^*|_{\pi_{Y^*}(Y^*)}||=||T||$ and $(T^*|_{\pi_{Y^*}(Y^*)})^*|_X=T.$
So, one can write $L(X,Y^{**})=L(Y^*,X^*)$ (in Banach sense).
\smallpagebreak

{\bf Lemma 2.}\
If $T\in L(X, Y)$ then

1)\, $\pi_Y T\in N_s(X,Y^{**}) \, \iff \, T^*\in N_s(Y^*, X^*);$

2)\, $T\in N_s(X,Y) \implies T^*\in N_s(Y^*, X^*).$
\smallpagebreak

{\bf Lemma 3.}\
If $E$ is a Banach space of type 2 (respectively, of cotype 2)
and of cotype $q_0$ (respectively, of type $q'_0)$
then $E$ has the $AP_s$ where $1/s=3/2-1/q_0.$
\smallpagebreak

The proofs of Lemma 3 can be found in [13], [14].          %RR!!!

\smallpagebreak

%%%%%%%%%%%%%%%%%%%%%%%%%%%%%%%%%%%  New for s<1 ...
%{\bf A theorem.}\
\smallpagebreak

\centerline{\bf \S3.\, A positive result}
\smallpagebreak

\proclaim {Theorem {\rm 1}}
Let $\, s\in (0,1],$ $ T\in L(X,Y)$
and either
$\, X^*\in \,AP_s\ $ or $\, Y^{***}\in \,AP_s.$
If $ T\in N_s(X, Y^{**}),$
then $T\in N_1(X,Y).$ In other words, under these conditions
from the $ s$-nuclearity of the conjugate operator $ T^*$ it follows
that the operator $ T$ belongs to the space $ N_1(X, Y)$ \, {\rm (that is nuclear)}.
\endproclaim

\demo{Proof}
Suppose there exists  an operator
    $ T\in L(X,Y)$
 such that
$ T\notin N_1(X,Y),$ but $ \pi_Y\,T\in N_s(X,Y^{**}).$
Since either $ X^*$ or $ Y^{**}$  has the $ AP_s,$
 $ N_s(X,Y^{**})=X^*\widehat\otimes_s Y^{**}.$
Therefore the operator $ \pi_Y\,T$ can be identified with the
tensor element
$ t\in X^*\widehat\otimes_s Y^{**}\subset X^*\widehat\otimes_1 Y^{**};$
in addition, by the choice of $ T,$ \
$  t\notin X^*\widehat\otimes_1 Y$ \ (the space
$  X^*\widehat\otimes_1 Y$
is considered as a closed subspace of the space
$  X^*\widehat\otimes_1 Y^{**}$\!).
Hence there is an operator
$ U\in L(Y^{**},X^{**})=\left( X^*\widehat\otimes_1 Y^{**}\right)^*$
with the properties that
$ \operatorname{trace}\, U\circ t=\operatorname{trace}\, \left(t^*\circ \left( U^*|_{X^*}\right) \right)=1$ and
$ \operatorname{trace}\, U\circ \pi_Y\circ z=0$ for each
$ z\in X^*\widehat\otimes_1 Y.$
From the last it  follows that, in particular, $ U\pi_Y=0$ and
$ \pi_Y^*\,U^*|_{X^*}=0.$
In fact, if
$ x'\in X^*$ and $ y\in Y,$ then
$$ <U\pi_Y\,y,x'> = <y, \pi_Y^*\,U^*|_{X^*}x'>
       = \operatorname{trace}\, \,U\circ (x'\otimes \pi_Y(y))=0.
$$
Evidently, the tensor element
$ U\circ t\in X^*\widehat\otimes_s X^{**}$  induces the operator
$ U\pi_Y T,$ which is equal identically to zero.

If $ X^*\in AP_s$ then
 $ X^*\widehat\otimes_s X^{**}= N_s(X,X^{**})$
and, therefore, this tensor element is zero what  contradicts
to the equality
$ \operatorname{trace}\,\, U\circ t=1.$

Let now $ Y^{***}\in AP_s.$ In this case
$$ V:= \left( U^*|_{X^*}\right)\circ T^*\circ \pi_Y^*: \
       Y^{***}\to Y^* \to X^*\to Y^{***}
$$
uniquely determines a tensor element
$ t_0$ from the $s$-projective tensor product
$ Y^{****}\widehat{\otimes}_s Y^{***}.$
Let us take any representation $ t=\sum x'_n\otimes y''_n$ for $ t$
as an element of the space $ X^*\widehat{\otimes}_s Y^{**}.$
Denoting for the brevity the operator $ U^*|_{X^*}$ by $ U_*,$
we obtain:
$$\multline
   Vy'''=U_*\, \left( T^*\pi_Y^*\,y'''\right) =
  %  U_*\, \left( (T^*\pi_Y^*\pi_{Y^*})\,\pi_Y^*\, y'''\right) =
   % U_*\, \left( (\pi_Y T)^*\,\pi_{Y^*})\,\pi_Y^*\, y'''\right) = \\ =
    U_*\, \left( (\sum y''_n\otimes x'_n) \,\pi_Y^*\, y'''\right)= \\
    =U_*\, \left( \sum <y''_n, \pi_Y^*\, y'''> \,x'_n \right) = \\ =
      \sum <\pi_Y^{**}y''_n,  y'''> \,U_* x'_n.
  \endmultline
$$

So, the operator $ V$ (or the element $ t_0$) has in the space
$ Y^{****}\widehat\otimes_s Y^{***}$ the representation
$$ V= \sum \pi_Y^{**}(y''_n)\otimes U_* (x'_n).
$$
Therefore,
$$  \operatorname{trace}\, t_0=\operatorname{trace}\, V= \sum <\pi_Y^{**}(y''_n), U_* (x'_n)> =
        \sum <y''_n, \pi_Y^*\,U_* x'_n> =  \sum 0=0
$$
(since $\pi_Y^*\,U_*=0;$ see above).

On the other hand,
$$  Vy'''= U_* \left( \pi_Y T\right)^* y'''= U_*\circ t^* (y''')=
     U_*\, \left( \sum <y''_n, y'''> \, x'_n\right)=
      \sum <y''_n, y'''> \, U_* x'_n,
$$
whence $ V=\sum y''_n\otimes U_*(x'_n).$   Therefore
$$ \operatorname{trace}\, t_0=\operatorname{trace}\, V= \sum <y''_n, U_* x'_n> = \sum <Uy''_n, x'_n>
= \operatorname{trace}\, U\circ t=1.$$
The obtained contradiction completes the proof of the theorem.
 $\quad\blacksquare$

\enddemo
%%%%%%%%%%%%%%%%%%%%%%%%%%%%%%%%%%%%%%%%%%%%%%%  End Th1 for s<1...

%%%%%%%%%%%%%%%%%%%%%%%%%%%%%%%%%%%%%%%%%%%%%%%%%%%%%
\vskip1.1cm
%{\bf Examples.}\

\smallpagebreak

\centerline{\bf \S4.\, Examples}
\smallpagebreak

We need two examples to show that all conditions, imposed on $X$ and $Y$ in Theorem 1,
are essential. Ideas of such examples are taken from the Reinov's work [10] (not translated,
as we know, from Russian).

\smallpagebreak

%\vskip0.1cm

{\bf Example 1.}\
Let $r\in(2/3,1], q\in[2,\infty), 1/r=3/2-1/q.$
There exist a separable reflexive Banach space $Y_0$ and a tensor element
$w\in Y_0^*\widehat\otimes_r Y_0$ so that
$w\neq0, \tilde w=0,$ the space $Y_0$ (as well as $Y_0^*)$
has the $AP_s$ for every $s<r$
(but, evidently, does not have the $AP_r).$
Moreover, $Y_0$ is of type 2 and of cotype $q_0$ for any $q_0>q.$
\vskip0.1cm

   %\demo{Proof}

{\it Proof}.\
%%%%%%%
We will use a variant of Per Enflo's [3] example of a Banach space without approximation property
given in the book [9, 10.4.5]. Namely, it follows from the constructions in [9] that           %RR!!!
there exist a Banach space $X$ and a tensor element $z\in X^*\widehat\otimes X$ so that
$\operatorname{trace}\, z\neq0,$ the operator $\tilde z,$ generated by $z,$ is identically zero and $z$ can be represented
in the following form:
$$(1) \quad  z=\sum_{N=1}^{\infty} \sum_{n=1}^{3\cdot 2^N}N^{1/2}2^{-3N/2}
  x'_{nN}\otimes x_{nN},
$$
where the sequences $(x'_{nN})$ and $ (x_{nN})$ are norm bounded by 1
(see [2] or [9, 10.4.5]).              %RR!!!

Fix $r\in(2/3,1]$ and
put  $1/q=3/2-1/r$ (thus, $q\in [2,\infty)$).
Let  $\{\varepsilon_N\}_1^{\infty}$ be a number sequence such that
 $\sum N^{-1-\varepsilon_N}<+\infty;$
$\gamma_N=2+3\varepsilon_N/2;$ $q_N$ is a number such that
$1/q-1/q_N=N^{-1}\log_2 N^{\gamma_N}$
(therefore, $q_N>q$).
Set $Y=\left(\sum_N l_{q_N}^{3\cdot 2^N}\right)_{l_q}.$

%%%%%%%%%

Denote by $e_{nN}$ (and $e'_{nN}$) the unit orths in $Y$ and in $Y^*$ respectively
($N=1,2,\dots; n=1,2,\dots, 3\cdot 2^N$), and put
$$
   z_1=\sum_N\sum_n 2^{-N/r}\,N^{-(1+\varepsilon_N)/r}\, e'_{nN}\otimes x_{nN};
   $$
   $$
   T=\sum_N\sum_n 2^{-(3/2-1/r)N}\, N^{(1/2+1/r+\varepsilon_N/r)}\, x'_{nN}\otimes e_{nN}.
$$

Let us show that $z_1\in Y^*\widehat\otimes_r X,$ $T\in \L(X,Y)$
(then, evidently, $z=z_1\circ T$\footnote{
  If $w\in Y^*\widehat\otimes X$ and $U\in L(X,Y)$ then $w\circ U$ denotes the image of $w$
  in $X^*\widehat\otimes X$ under the map $U^*\otimes \operatorname{id}_X.$
  So, if $w=\sum \varphi'_m\otimes \psi_m$ is a representation of $w$ in $Y^*\widehat\otimes X$ then
  $w\circ U=\sum U^*\varphi'_m\otimes \psi_m$ is a representation of $w\circ U$ in $X^*\widehat\otimes X.$
  Then $\operatorname{trace}\, w\circ U= \sum \langle U^*\varphi'_m,  \psi_m\rangle= \sum \langle \varphi'_m, U\psi_m\rangle= \operatorname{trace}\, U\circ w,$
  where $U\circ w$ is the image of $w$ in $Y^*\widehat\otimes Y$ under the map $\operatorname{id}_{Y^*}\otimes U.$
  }).
The first inclusion is evident (because of the choice of $\varepsilon_N).$

If $||x||\leqslant1$ then
$$
||Tx||\leqslant 3\left(\sum_N \left(N^{1/2+1/r+\varepsilon_N/r}/2^{(3/2-1/r-1/q_N)N}\right)^q\right)^{1/q}.
$$
Since $3/2-1/r-1/q_N=1/q-1/q_N=N^{-1}\log_2 N^{\gamma_N},$
we get from the last inequality:
$$
||T||\leqslant 3\left(\sum N^{(1/2+1/r+\varepsilon_N/r-\gamma_N)q}\right)^{1/q}=
       3 \left(\sum N^{-1-\varepsilon_N}\right)^{1/q}<\infty.
$$

Hence,
$$
 z=z_1\circ T,\ \ \tilde{z}:\ X\overset{T}\to\to Y\overset{\tilde{z}_1}\to\to X.
$$

Now, let $Y_0:=\overline{T(X)}\subset Y, $\, $T_0: X\to Y_0$ be induced by $T$ and
$z_0:= z_1\circ j$ where $j: Y_0\hookrightarrow Y$ is the natural embedding.
Then $T_0\in L(X,Y_0), z_0\in Y_0^*\widehat\otimes_r X,$
$z=z_1\circ T=z_0\circ T_0,$ $\operatorname{trace}\, z_0\circ T_0\neq0$
(so, $z_0\neq0)$  and $\tilde z_0=0.$

Write $z_0$ as $z_0=\sum_{m=1}^\infty f'_m\otimes f_m,$ where
$(f'_m)\subset Y_0^*, (f_m)\subset X$ and $\sum_m ||f'_m||^r ||f_m||^r<\infty.$
We get:
$$
\operatorname{trace}\, z_0\circ T_0 = \sum \langle T_0^*f'_m, f_m\rangle = \sum \langle f'_m, T_0f_m\rangle= \operatorname{trace}\, T_0\circ z_0.
$$
Therefore, $w:= T_0\circ z_0\in Y_0^*\widehat\otimes_r Y_0,$ $\operatorname{trace}\, w\neq0$ and $\tilde w=0.$

Since the space $Y$ is of type 2 and of cotype $q_0$ for every $q_0>q$ (and $Y^*$ is
of cotype 2 and of type $q'_0$),
the space $Y_0$ (respectively, $Y_0^*)$ has the $AP_s,$
where $1/s=3/2-1/q_0,$ for every $s<r$ (Lemma 3).

  %\enddemo

%%%%%%%%%%%%%%

 \vskip0.1cm
%%%%%%%%%%%%%%%%%%%%%%%%%%%%%%%%%%%%%

{\it Remark}:\
We have a nice "by-product consequence" of Example 1.
For $q=2$ (that is, $r=1)$, the space $Y_0$ is a subspace of the space of type
$\left(\sum_j l_{p_j}^{k_j}\right)_{l_2}$ with $p_j\searrow 2$ and $k_j\nearrow \infty.$ 
Every such space is an asymptotically Hilbertian space (for definitions and some discussion,
see [1]). So, we got:
\smallpagebreak

{\bf Corollary.}\
There exists an asymptotically Hilbertian space without the Grothendieck approximation property.
 \smallpagebreak
 
 First example of such a space was constructed (by O. Reinov) in 1982 [10],          %R!!!
 where A. Szankowski's results were used
 (let us note that in that time there was not yet such  notion as "asymptotically Hilbertian space"). 
 Later,  in 2000, by applying Per Enflo's example in a version of  A.M. Davie [2],        %R!!!
 P.~G. ~Casazza,   C.~L. ~Garc\'{\i}a and   W.~B. ~Johnson [1]                       %R!!!
 gave another example of an asymptotically Hilbertian space  which fails the approximation property.
 We here, not being searching for an example of such a space, have got it (accidentally)
 by using the construction from [9].           %R!!!

% Banach spaces without approximation property
% Reinov O.I., 
% Funktsional'nyi Analiz i Ego Prilozheniya, Vol. 16, No. 4, pp. 84-85, October-December, 1982.
%     Functional Analysis and Its Applications
%October–December, 1982, Volume 16, Issue 4, pp 315-317
%Banach spaces without approximation property
%O. I. Reinov

 % An example of an asymptotically Hilbertian space  which fails the approximation property
%[10] PG. Casazza, C.L. Garcia and W.B. Johnson, An example of an asymptotically hilbertian space which fails
%the approximation property, Preprint.
   % P.~G. ~Casazza   C.~L. ~Garc\'{\i}a          W.~B. ~Johnson
% arXiv:math/0006134v2 [math.FA] 20 Sep 2000
                          % Proceedings of the American Mathematical Society, ISSN 0002-9939, Vol. 129, N? 10, 2001 , pags. 3017-3024
%                          EXAMPL 2

\vskip0.1cm

{\bf Example 2.}\
Let $r\in [2/3,1), q\in(2,\infty], 1/r=3/2-1/q.$
There exist a subspace $Y_q$ of the space $l_q$ and a tensor element
$w_q\in Y_q^*\widehat\otimes_1 Y_q$ so that
$w_q\in Y_q^*\widehat\otimes_s Y_q$ for each $s>r,$
$w_q\neq0, \tilde w_q=0$  and the space $Y_q$ (as well as $Y_q^*)$
has the $AP_r$
(but, evidently, does not have the $AP_s$ if $1\ge s>r).$
Clearly, $Y_q$ is of type 2 and of cotype $q$ for $q<\infty.$
\vskip0.1cm

     %\demo{Proof}
{\it Proof}.\
We are going to follow the way  indicated in the proof of the assertion from Example 1.
Let $X$ and $z$ be as in that proof, so that $z$ has the form (1).
Fix $r\in [2/3,1).$ Now $1/q=3/2-1/r,$ and we put $\varepsilon_N=0$ and $\gamma_N=2$
for all $N;$ all $q_n$'s are equal to $q.$ Let us fix also an $\alpha=\alpha(q)>0$ (to be defined later).
Consider the space $Y:= \left(\sum_N l_{q}^{3\cdot 2^N}\right)_{l_q}$ (in the case $q=\infty$
"$l_q$" means "$c_0$").

Denote by $e_{nN}$ (and $e'_{nN}$) the unit vectors  in $Y$ and in  $Y^*$ respectively
($N=1,2,\dots; n=1,2,\dots, 3\cdot 2^N$), and set this time

$$
   z_1=\sum_N\sum_n 2^{-N/r}\,N^{\alpha}\, e'_{nN}\otimes x_{nN};
   $$
   $$
   T=\sum_N\sum_n 2^{-(3/2-1/r)N}\, N^{1/2-\alpha}\, x'_{nN}\otimes e_{nN}.
$$

Then $z_1\in Y^*\widehat\otimes_{s} X$ for every $s>r.$
Indeed, if $s\in (r,1]$ then
$$
 \sum_N\sum_n [2^{-N/r}\,N^{\alpha}]^{s}=  \sum_N 3\cdot 2^N\cdot 2^{-Ns/r}\,N^{\alpha\, s}=
 \sum_N 3\cdot 2^{-\varepsilon_0 N}\, N^{\alpha s}<\infty,
$$
where $\varepsilon_0=s/r-1>0.$

Show that  $T\in \L(X,Y)$ (clearly, then $z=z_1\circ T$).
Indeed,
if $||x||\leqslant1$ then, for $q<\infty:$
$$
 ||Tx||_{l_q}^q\le
 \sum_N\sum_n \left(2^{-N/q}\, N^{1/2-\alpha}\,  |\langle x'_{nN},x\rangle|\right)^q\le
  \sum_N 3\cdot 2^{N-(N/q)q}\, N^{(1/2-\alpha)q} = \sum_N 3\cdot 1\cdot N^{\alpha_0},
$$
where $\alpha_0=(1/2-\alpha)q.$
Now, take $\alpha>0$ such that  $\alpha_0=-2.$
For $q=\infty$ take $\alpha=\alpha(\infty)=1.$

Therefore,
$$
 z=z_1\circ T,\ \ \tilde{z}:\ X\overset{T}\to\to Y\overset{\tilde{z}_1}\to\to X;
$$

As in the case of the previous proof (in Example 1),
 let $Y_q:=\overline{T(X)}\subset Y, T_q: X\to Y_q$ be induced by $T$ and
$z_q:= z_1\circ j$ where $j: Y_q\hookrightarrow Y$ is the natural embedding.
Then $T_q\in L(X,Y_q), z_q\in Y_q^*\widehat\otimes_s X$ for all $s>r,$
$z=z_1\circ T=z_q\circ T_q,$ $\operatorname{trace}\, z_q\circ T_q\neq0$
(so, $z_q\neq0)$  and $\tilde z_q=0.$

Write $z_q$ as $z_q=\sum_{m=1}^\infty f'_m\otimes f_m,$ where
$(f'_m)\subset Y_q^*, (f_m)\subset X$ and $\sum_m ||f'_m||\, ||f_m||<\infty.$
We get:
$$
\operatorname{trace}\, z_q\circ T_q = \sum \langle T_q^*f'_m, f_m\rangle = \sum \langle f'_m, T_qf_m\rangle= 
\operatorname{trace}\, T_q\circ z_q.
$$
Therefore, $w_q:= T_q\circ z_q\in Y_q^*\widehat\otimes_s Y_q$ for every $s>r,$ 
$\operatorname{trace}\, w_q\neq0$ and $\tilde w_q=0.$
Finally, Lemma 3 says that the space $Y_q$ has the $AP_r$ if $q<\infty.$
If $q=\infty$ then, as we know, any Banach space has the $AP_{2/3}.$
        %\enddemo
\smallpagebreak

%%%%%%% end EXample 2

{\it Remark}:\
The space $Y_\infty$ from Example 2 not only does not have the $AP_s$
for any $s\in (2/3, 1],$ but also does not have the $AP_p$ (in the sense of the paper [12])    %R!!!
for all $p\in[1,2)$ (this follows from some facts proved in [13]).          %R!!!

\smallpagebreak

\vskip0.8cm

\centerline{\bf  \S4.\, Applications of Examples}
\vskip0.5cm

Next two theorems show that the conditions "$X^*$ has the $AP_s$" and
"$Y^{***}$ has the $AP_s$" are essential in Theorem 1 and can not be replaces by weaker conditions
"$X$ has the $AP_s$" (even by "$X^*$ has the $AP_1$") or "$Y^{**}$ has the $AP_s$";
moreover, even "both $X$ and $Y^{**}$ have the $AP_1$" is not enough for the conclusion
of Theorem 1 to be valid.
\smallpagebreak

{\bf Theorem 2.}\
Let $r\in (2/3,1], q\in [2,\infty), 1/r=3/2-1/q.$
There exist a Banach space $Z_0$ and an operator $T\in L(Z_0^{**}, Z_0)$ so that

(1)\,
$Z_0^{**}$  has a Schauder basis (so, has the MAP);

(2)\,
all the duals of $Z_0$ are separable;

(3)\,
$Z_0^{***}$ has the $AP_s$ for every $s\in (0,r);$

(4)\,
$\pi_{Z_0} T\in N_r(Z_0^{**}, Z_0^{**});$

(5)\,
$T\notin N_1(Z_0^{**}, Z_0);$

(6)\,
$Z_0^{***}$ does not have the $AP_r.$

\smallpagebreak

\demo{\it Proof}
Let us fix $r\in (2/3,1], q\in [2,\infty), 1/r=3/2-1/q$\, and take the pair $ (Y_0, w)$
from Example 1.
Let $ Z_0$ be a separable space such that $ Z_0^{**}$ has a basis
and there exists a linear homomorphism $ \varphi$ from $ Z_0^{**}$
onto $ Y_0$ with the kernel
$ Z_0\subset Z_0^{**}$ so that the subspace $ \varphi^*(Y_0^*)$ in complemented
in $ Z_0^{***}$ and, moreover,
$Z_0^{***}\cong \varphi^*(Y_0^*)\oplus Z_0^*$
(see [7, Proof of Corollary 1]). Lift the tensor element
$ w,$ lying in $ Y_0^*\widehat\otimes_r Y_0,$ up to an element \footnote{
   If $ w=\sum_{k=1}^\infty \,y'_k\otimes\,y_k $ is any representation of $ w$
in  $ Y_0^*\widehat\otimes_r Y_0,$  then we take $ \{ z''_n\}\subset Z_0^{**}$
in such a way that
the last sequence is absolutely $ r$-summing and
$ \varphi(z''_n)=y_n$ for every $ n.$
   }
$ w_0\in Y_0^*\widehat\otimes_r Z_0^{**},$
so that $ \varphi\circ w_0=w,$ and set $ T:= w_0\circ \varphi.$
Since $ \operatorname{trace}\, w_0\circ \varphi=\operatorname{trace}\, \varphi\circ w_0=\operatorname{trace}\, w=1$ and $ Z_0^{**}$ has the
AP, then $ \widetilde {w_0}=w_0\neq 0.$
Besides, the operator $ \widetilde{\varphi\circ w_0}:Y_0\to Z_0^{**}\to Y_0,$
associated with the tensor $ \varphi\circ w_0,$ is equal to zero. Therefore
$ w_0(Y_0)\subset \operatorname{ Ker}\varphi= Z_0\subset Z_0^{**},$
that is the operator
$ w_0$ is acted from $ Y_0$ into $ Z_0.$

Since the subspace $ \varphi^*(Y_0^*)$ is complemented in $ Z_0^{***},$
then $ w_0\circ\varphi\in Z_0^{***}\widehat\otimes_r Z_0= N_r(Z_0^{**},Z_0)$ iff
$ w_0\in Y_0^{*}\widehat\otimes_r Z_0= N_r(Y_0,Z_0).$

If $ w_0\in  N_r(Y_0,Z_0),$
then, for its arbitrary (nonzero!) $N_r$-representation
of the form $ w_0=\sum y'_n\otimes z_n,$ the composition $ \varphi\circ w_0$ is
a zero tensor element in
$ Y_0^*\widehat\otimes_r Y_0;$ but this composition
represents the element $ w, $ which, by its choice, can not be zero.
Thus, $ w_0\notin  N_r(Y_0,Z_0)$ and, thereby,
 $ w_0\circ\varphi\notin Z_0^{***}\widehat\otimes_r Z_0= N_r(Z_0^{**},Z_0).$
On the other hand, certainly,
 $ w_0\circ\varphi\in Z_0^{***}\widehat\otimes_r Z_0^{**}= N_r(Z_0^{**},Z_0^{**}).$

 Finally, since $Z_0^{***}\cong \varphi^*(Y_0^*)\oplus Z_0^*,$ one has that
 the space $Z_0^{***}$ has the $AP_s$ for every $s\in (0,r).$
 $\quad\blacksquare$
\enddemo

%%   end TH 2 PROOF

%{\it Remark}.  [ABOUT N-P]

   %%%%%%%%%%%%%%%  TH 3
\smallpagebreak

{\bf Theorem 3.}\
Let $r\in [2/3,1), q\in (2,\infty], 1/r=3/2-1/q.$
There exist a Banach space $Z_q$ and an operator $T\in L(Z_q^{**}, Z_q)$ so that

(1)\,
$Z_q^{**}$  has a Schauder basis (so, has the MAP);

(2)\, if $q<\infty$ then
all the duals of $Z_q$ are separable;

(3)\,
$Z_q^{***}$ has the $AP_r;$

(4)\,
$\pi_{Z_q} T\in N_s(Z_q^{**}, Z_q^{**})$ for every $s\in (r,1];$

(5)\,
$T\notin N_1(Z_q^{**}, Z_q);$

(6)\,
$Z_q^{***}$ does not have the $AP_s$ for any $s\in (r,1];$

\smallpagebreak

\demo{\it Proof}
Let us fix $r\in [2/3,1), q\in (2,\infty], 1/r=3/2-1/q$\, and take the pair $ (Y_q, w_q)$
from Example 2.
Let $ Z_q$ be a separable space such that $ Z_q^{**}$ has a basis
and there exists a linear homomorphism $ \varphi$ from $ Z_q^{**}$
onto $ Y_q$ with the kernel
$ Z_q\subset Z_q^{**}$ so that the subspace $ \varphi^*(Y_q^*)$ in complemented
in $ Z_q^{***}$ and, moreover,
$Z_q^{***}\cong \varphi^*(Y_q^*)\oplus Z_q^*$
(as in the proof of Theorem 2, see [7, Proof of Corollary 1]).

Construct $ w_0\in Y_q^*\widehat\otimes_1 Z_q^{**}$ (following the way of the proof
of Theorem 2)
so that $ w_0\in Y_q^*\widehat\otimes_s Z_q^{**}$ for every $s\in (r,1]$ and
$ \varphi\circ w_0=w_q$
(it is possible to apply a "simultaneous lifting" procedure --- see Footnote 2, --- 
since $w$ has a form from the proof of the assertion of Example 2). 
Set $ T:= w_0\circ \varphi.$
From this point, the proof repeats the arguments of the proof of Theorem 2,
and we have to mention only:
since $Z_q^{***}\cong \varphi^*(Y_q^*)\oplus Z_q^*,$ one has that
 the space $Z_q^{***}$ has the $AP_r.$

\enddemo

%%   end TH 3 PROOF

%%

%{\it Remark}.  [ABOUT N-P]

%\newpage

\bigpagebreak

\centerline{REFERENCES}
\smallpagebreak

\ref \no 1  \by  Casazza  P.~G., Garc\'{\i}a  C.~L.,  Johnson W.~B.,\pages 3017-3024
\paper An example of an asymptotically Hilbertian space  which fails the approximation property
\yr 2001 \vol 129 \issue 10
\jour Proceedings of the American Mathematical Society
\endref

\ref \no 2\by Davie A.M. \pages  261--266
\paper  The approximation problem for Banach spaces
\yr 1973\vol 5
\jour Bull. London Math. Soc.
\endref

\ref \no 3  \by Enflo P.  \pages  309--317
\paper A counterexample to the approximation property in Banach spaces
\yr 1973 \vol 130 \issue
\jour  Acta Math.
\finalinfo  %$MR 53 \# 6288.$
\endref

\ref \no 4\by Figiel T., Johnson W.B.\pages 197--200
\paper  The approximation property does not imply  the bounded
   approximation property
\yr 1973\vol 41
\jour   Proc. Amer. Math. Soc.
\endref

\ref \no5\by Grothendieck A. \pages 196 + 140
\paper  Produits tensoriels topologiques et espases nucl\'eaires
\yr 1955\vol  16
\jour  Mem. Amer. Math. Soc.
\endref

\ref \no 6\by  Hinrichs A.,  Pietsch A.\pages  232--261
\paper  $p$-nuclear operators in the sense of Grothendieck
\yr 2010\vol 283 \issue 2
\jour Math. Nachr. 
\endref

\ref \no 7\by Lindenstrauss J.\pages  279--284
\paper  On James' paper ``Separable Conjugate Spaces"
\yr 1971\vol 9
\jour Israel J. Math.
\endref

\ref \no 8  \by Oja E., Reinov O.I.   \pages  121--122
\paper  Un contre-exemple \`a une affirmation de A.Grothendieck
\yr  1987 \vol  305
\jour   C. R. Acad. Sc. Paris. --- Serie I
\endref

\ref \no 9 \by Pietsch  A.\pages 536 p
 \paper  Operator ideals
 \yr 1978\vol
 \jour North-Holland
 \endref

   \ref \no 10  \by Reinov  O. I. \pages 315-317
\paper Banach spaces without approximation property
\yr 1982 \vol 16 \issue 4
\jour Functional Analysis and Its Applications
\endref

\ref \no 11  \by Reinov O.I. \pages  115-116
\paper  A simple proof of two theorems of A. Grothendieck
\yr 1983 \vol 7
\jour Vestn. Leningr. Univ.
\endref

\ref \no 12  \by Reinov  O. I. \pages   125-134
\paper   Approximation properties of order p and the existence of
  non-p-nuclear operators with p-nuclear second adjoints
\yr 1982 \vol  109
\jour    Math. Nachr.
\endref

\ref \no 13  \by Reinov O.I.\pages 145-165
\paper   Disappearing tensor elements in the scale of p-nuclear operators
\yr 1983\vol 1
\jour  in: Theory of Operators and Theory of Functions, Leningrad , LGU
\endref

\ref \no 14  \by Reinov  O. I. \pages 2243-2250
\paper Approximation properties $ \operatorname{AP_s}$ and $p$-nuclear
operators {\rm(}the case  $ 0<s\le1)$
\yr 2003 \vol 115 \issue 2
\jour Journal of Mathematical Sciences
\endref

\enddocument

                        \comment

% Banach spaces without approximation property
% Reinov O.I., 
% Funktsional'nyi Analiz i Ego Prilozheniya, Vol. 16, No. 4, pp. 84-85, October-December, 1982.
%     Functional Analysis and Its Applications
%October–December, 1982, Volume 16, Issue 4, pp 315-317
%Banach spaces without approximation property
%O. I. Reinov

 % An example of an asymptotically Hilbertian space  which fails the approximation property
%[10] PG. Casazza, C.L. Garcia and W.B. Johnson, An example of an asymptotically hilbertian space which fails
%the approximation property, Preprint.
   % P.~G. ~Casazza   C.~L. ~Garc\'{\i}a          W.~B. ~Johnson   
% arXiv:math/0006134v2 [math.FA] 20 Sep 2000 
                          % Proceedings of the American Mathematical Society, ISSN 0002-9939, Vol. 129, N? 10, 2001 , pags. 3017-3024

  %%%%%%%%%%%%%%%%%%%%%%%%%%%%%%%%%%%%%%%%
O. I. Reinov
Approximation properties $AP_s$ and $p$-nuclear operators
(the case $0<s\le 1)$
Journal of Mathematical Sciences, Vol. 115, No. 2, 2003
2243--2250

A. M. Davie, \The approximation problem for Banach spaces," Bull. London Math. Soc., 5, 261{266 (1973).
2. T. Figiel and W. B. Johnson, \The approximation property does not imply the bounded approximation
property," Proc. Am. Math. Soc., 41, 197{200 (1973).
3. A. Grothendieck, \Produits tensoriels topologiques et espaces nucleaires," Mem. Am. Math. Soc., 16, 196{140
(1955).
4. J. Lindenstrauss, \On James' paper `Separable Conjugate Spaces'," Israel J. Math., 9, 279{284 (1971).

6. E. Oja and O. I. Reinov, \Un contre-exemple a une armation de A. Grothendieck," C. R. Acad. Sci. Paris.
Serie I, 305, 121{122 (1987).
7 A. Pietsch, Operator Ideals, North-Holland (1980).
8. O. I. Reinov, \Approximation properties of order p and the existence of non-p-nuclear operators with p-nuclear
second adjoints," Math. Nachr., 109, 125{134 (1982).
9 O. I. Reinov, \A simple proof of two theorems of A. Grothendieck," Vestn. Leningr. Univ., 7, 115{116 (1983).
10. O. I. Reinov, \Disappearing tensor elements in the scale of p-nuclear operators," in: Theory of Operators and
Theory of Functions, Leningrad (1983), pp. 145{165.

%%%%%%
\ref \no 3  \by O. I. Reinov  \pages 277-291\nofrills
\paper
Approximation properties $ \operatorname{AP_s}$ and $p$-nuclear
operators {\rm(}the case when $ 0<s\le1)$
\yr 2000 \vol 270
\jour Zapiski nauchn. sem. POMI
\finalinfo (in Russian)
\endref

\ref \no 4\by J. Lindenstrauss \pages  279-284
\paper  On James' paper "Separable Conjugate Spaces"
\yr 1971\vol 9
\jour Israel J. Math.
\endref

\ref \no 1\by  Hinrichs A.,  Pietsch A.\pages  232-261
\paper  $pp$-nuclear operators in the sense of Grothendieck
\yr 2010\vol 283
\jour Math. Nachr. 
\endref

\ref \no 2\by Figiel T., Johnson W.B.\pages 197-200
\paper  The approximation property does not imply  the bounded
   approximation property
\yr 1973\vol 41
\jour   Proc. Amer. Math. Soc.
\endref

\ref \no3\by Grothendieck A. \pages 196-140
\paper  Produits tensoriels topologiques et espases nucl\'eaires
\yr 1955\vol  16
\jour  Mem. Amer. Math. Soc.
\endref

\ref \no 4\by Lindenstrauss J.\pages  279-284
\paper  On James' paper ``Separable Conjugate Spaces"
\yr 1971\vol 9
\jour Israel J. Math.
\endref

\ref \no 5\by ia¬a?«o u.i., ?a?a??¬?? e.a.\pages   122-144
\paper  ??aia? ?N¬oN?µ?a???a? »«??«?a ? i???¬?N ¬ ?N? ?o«???oa ?N¬«?«???
  »?«???a???o «»N?a?«?«o
\yr 1983\vol 1
\jour o ¬?. ``?N«??? «»N?a?«?«o ? ?N«??? ???¬µ??". i.:iao
\endref

\ref \no 6  \by Oja E., Reinov O.I.   \pages  121-122
\paper  Un contre-exemple \`a une affirmation de A.Grothendieck
\yr  1987 \vol  305
\jour   C. R. Acad. Sc. Paris. --- Serie I
\endref

\ref \no 7 \by A?? C.\pages 536 ?
 \paper A»N?a?«???N ?nNa??
 \yr 1982\vol
 \jour     i«?¬oa: i??
 \endref

\ref \no 8  \by Reinov O.I.\pages   125-134
\paper   Approximation properties of order p and the existence of
  non-p-nuclear operators with p-nuclear second adjoints
\yr 1982 \vol  109
\jour    Math. Nachr.
\endref

\ref \no 9  \by Reinov O.I. \pages  115-116
\paper  A simple proof of two theorems of A. Grothendieck
\yr 1983 \vol 7
\jour Vestn. Leningr. Univ.
\endref

\ref \no 10  \by Reinov O.I.\pages 145-165
\paper   Disappearing tensor elements in the scale of p-nuclear operators
\yr 1983\vol 1
\jour  in: Theory of Operators and Theory of Functions, Leningrad , LGU
\endref

9 O. I. Reinov, \A simple proof of two theorems of A. Grothendieck," Vestn. Leningr. Univ., 7, 115{116 (1983).
10. O. I. Reinov, \Disappearing tensor elements in the scale of p-nuclear operators," in: Theory of Operators and
Theory of Functions, Leningrad (1983), pp. 145{165.

Aicke Hinrichs and Albrecht Pietsch
p -nuclear operators in the sense of Grothendieck
Math. Nachr. 283, No. 2, 232 – 261 (2010)
            \endcomment

\enddocument